\documentclass[oneside,english,reqno]{amsart}
\usepackage{palatino}
\usepackage[T1]{fontenc}
\usepackage[latin1]{inputenc}
\usepackage{amssymb}
\IfFileExists{url.sty}{\usepackage{url}}
                      {\newcommand{\url}{\texttt}}

\makeatletter

 \theoremstyle{plain}
 \theoremstyle{definition}
 \newtheorem*{defn*}{Definition}
 \theoremstyle{plain}    
 \newtheorem*{thm*}{Theorem} 
 \theoremstyle{plain}    
 \newtheorem{lem}{Lemma} 

\ifx\hyperlink\undefined
\newcommand{\hrlb}{ }
\else
\newcommand{\hrlb}{\\}
\fi

\usepackage{babel}
\makeatother
\begin{document}

\title{Excited random walk in two dimensions has linear speed}

\author{Gady Kozma}

\thanks{This material is based upon work supported by the National Science
Foundation under agreement DMS-0111298. Any opinions, findings and
conclusions or recommendations expressed in this material are those
of the authors and do not necessarily reflect the views of the National
Science Foundations}

\address{Institute for Advanced Study, 1 Einstein dr., Princeton NJ 08540
USA}

\email{gady@ias.edu}

\maketitle
Imagine a cookie placed on every vertex of an infinite $d$-dimensional
grid. A random walker on this grid behaves as follows upon encountering
a cookie: he consumes it and then performs a random step with an $\epsilon$-drift
to the right, namely the probability to make a right move is $\frac{1}{2d}+\epsilon$,
the probability for a left move is $\frac{1}{2d}-\epsilon$ and the
probability for all other moves is $\frac{1}{2d}$. When encountering
a site already visited (so no cookie), he performs a simple random
walk. This process was coined by Benjamini and Wilson \cite{BW03}
{}``excited random walk'' (though the name {}``brownie walk''
might describe it better). Since then a number of papers were devoted
to this process. See \cite{V03,K,Z,ABK}, \cite{PW97,D99} for a one
dimensional continuous version and \cite{AR} for some simulation
results.

\cite{BW03} mostly discusses the geometric case $d>1$. They prove
that in dimension $\geq4$ the walk has linear speed, namely \begin{equation}
\varliminf_{n\to\infty}\frac{E(n)_{1}}{n}>0\label{eq:speed}\end{equation}
where $E(n)$ is the position of the walk at time $n$ and $E(n)_{1}$
is its first (left-right) coordinate. In dimension $2$ they prove
transience, in fact they prove that $E(n)_{1}>cn^{3/4}\log^{-5/4}n$
for $n$ sufficiently large almost surely, and ask what is the correct
speed. \cite{K} extended (\ref{eq:speed}) to the three dimensional
case. The purpose of this paper is to show the same for two dimensions
(in one dimension this is not true, though the multiple cookies case
discussed by Zerner \cite{Z} is still open).

Benjamini and Wilson's proof of transience in two dimensions will
play a crucial role in the current paper so let's describe it briefly.
They coupled excited random walk to a simple random walk $R$ in the
natural way: when the ERW encounters a cookie they walk {}``as close
as possible'', that is with probability $2\epsilon$ the ERW walks
to the right and the SRW walks to the left, while with probability
$1-2\epsilon$ they perform the same step. If no cookie, they just
perform the same step. This means that (with this coupling) the ERW
is always to the right of the SRW and the distance is increasing.
This implies that when the SRW reaches a tan point, a point $(x,y)$
such that the random walk never visited before $(x+n,y)$ for any
$n=0,1,2,\dotsc$ then the ERW must be at a vertex with a cookie.
The name {}``tan point'' comes from placing the sun at right infinity,
and these points are the points when the SRW can tan since its past
path is not blocking the sun. Hence an estimate of the number of tan
points gives a lower bound for the number of cookies eaten. 

Thus we are left with a problem on simple random walk, i.e.~estimate
the probability that $R(n)$ is a tan point. By symmetry this is the
same as the probability that a random walk of length $n$ will avoid
hitting a half line. This, as is well known is $\approx n^{-1/4}$
--- see \cite{K87}. The same $\frac{1}{4}$ translates to the $\frac{3}{4}$
exponent in the speed estimate of Benjamini and Wilson with some logarithmic
corrections.

To get better than $n^{3/4}$ one needs to apply ERW's {}``self-correcting''
property. Very roughly, if an ERW of length $m$ goes to a distance
of $\mu$ for $\mu\gg\sqrt{m}$ then the next portion of length $m$
of the ERW should be quite independent from the previous portion,
and should also continue to a length of $\mu$. Unfortunately, one
cannot continue this forever and claim that the speed is $>\mu/m$
because there is always a small probability for a portion to fail
and then one needs a {}``fallback mechanism''. This is provided
by the Benjamini and Wilson argument (performed locally, see the definition
of a relative tan point below). Thus the proof is inductive, using
Benjamini and Wilson to both kickstart the induction and to provide
a fallback mechanism in each stage.

I wish to thank Itai Benjamini and Martin Zerner for enlightening
discussions of this problem.

\section{Preliminaries}

\begin{defn*}
Let $\epsilon<\frac{1}{4}$. Let $\mathcal{V}\subset\mathbb{Z}^{2}$
be some subset of vertices ($\mathcal{V}$ standing for {}``visited'')
and let $w\in\mathbb{Z}^{2}$. $\epsilon$-Excited random starting
from $(\mathcal{V},w)$ is a stochastic process $E(n)$ on $\mathbb{Z}^{2}$
such that $E(0)=w$ and such that if $E(n)=E(m)$ for some $m<n$
or $E(n)\in\mathcal{V}$ then $E(n+1)$ has probability $\frac{1}{4}+\epsilon$
to be $E(n)+(1,0)$, probability $\frac{1}{4}-\epsilon$ to be $E(n)+(-1,0)$
and probability $\frac{1}{4}$ to be either $E(n)+(0,\pm1)$. In the
other case, that $E(n)\not\in E[0,n-1]\cup\mathcal{V}$, distribute
$E(n+1)$ like a simple random walk.
\end{defn*}
Let us next state the result with an explicit bound on the probability:

\begin{thm*}
Let $E$ be an $\epsilon$-excited random walk starting from $(\emptyset,(0,0))$.
Then with probability $>1-Ce^{-c\log^{2}n}$ one has $E(n)_{1}>cn$.
The constants $C$ and $c$ may depend on $\epsilon$.
\end{thm*}
For a simple or excited walk $R$ and two times $i\leq j$ we shall
denote by $R[i,j]$ the set $\{ R(k):i\leq k\leq j\}\subset\mathbb{Z}^{2}$.
$\log x$ will always be a shortcut for $\max\{\log x,1\}$. $\left\lfloor x\right\rfloor $
and $\left\lceil x\right\rceil $ will denote, as usual, the largest
integer $\leq x$ and the smallest integer $\geq x$ respectively.
By $C$ and $c$ we shall denote constants depending only on $\epsilon$
whose precise value is unimportant as far as this paper is concerned,
and could change from formula to formula or even within the same formula.
$C$ will pertain to constants which are {}``big enough'' and $c$
to constants which are {}``small enough''. We will number a few
$C$ and $c$-s --- only those which we will reference later on. When
we say {}``$x$ is sufficiently large'' we mean {}``$x>C$ for
some $C$'' and in particular the bound may depend on $\epsilon$.

\section{Proof}

\begin{defn*}
\label{def:tanp}Let $R$ be a simple random walk and let $i<j$ be
two times. Then we say that SRW has a \textbf{tan point} at $j$ \textbf{relative}
to $i$ if the portion of the walk $R[i,j-1]$ does not intersect
the half line $R(j)+\mathbb{N}\times\{0\}$.
\end{defn*}
\begin{lem}
\label{lem:imj}Let $R$ be a simple random walk of length $n$, and
let $m<n$ be sufficiently large. Then with probability $\geq1-Ce^{-c\log^{2}n}$
the following holds: for every $0\leq i<n-m\log^{6}n$ the random
walk on $\left[i,i+\left\lfloor m\log^{6}n\right\rfloor \right]$
exhibits $\geq c_{1}m^{3/4}$ tan points relative to $i$ in some
interval $[j,j+m]$ for $i+m\leq j\leq i+m\log^{6}n-m$.
\end{lem}
\begin{proof}
Fix $i$. Denote by $B_{1}$ ($B$ standing for {}``bad'') the event
that \[
\max_{j\in[i,i+m]}|R(j)_{2}-R(i)_{2}|>{\textstyle \frac{1}{2}}\log^{2}n\sqrt{m}.\]
 It is well known that $\mathbb{P}(B_{1})\leq Ce^{-c\log^{4}n}$ ---
for each $j$ this follows from the Chernoff bound and summing over
$j$ only changes the constant in the exponent. Next denote by $B_{2}$
the event \[
\max\left\{ |R(j)_{2}-R(i)_{2}|:j\in\left[i,i+\left\lfloor m\log^{6}n\right\rfloor -m\right]\right\} <2\log^{2}n\sqrt{m}.\]
 Again, it is well known that $\mathbb{P}(B_{2})\leq Ce^{-c\log^{2}n}$
--- for any $x\in I$, \[
I:=\left[\left\lceil R(i)_{2}-2\log^{2}n\sqrt{m}\right\rceil ,\left\lfloor R(i)_{2}+2\log^{2}n\sqrt{m}\right\rfloor \right]\]
there is a probability $>c$ to exit $I$ in the next $\left\lfloor m\log^{4}n\right\rfloor $
steps (if $m$ is sufficiently large), and so the probability to exit
it in $\left\lfloor m\log^{6}n\right\rfloor -m$ steps is $\geq1-Ce^{-c\log^{2}n}$.

Let $\lambda$ be some parameter to be fixed later and assume $m>1/\lambda$.
Let $k\in\mathbb{Z}\setminus\{0\}$ be positive or negative and denote
by $L_{k}$, $S_{k}$ and $H_{k}$ the horizontal line, strip and
half strip at height $R(i)_{2}+k\left\lfloor \sqrt{\lambda m}\right\rfloor $
respectively. In a formula:\begin{align*}
L_{k} & :=\mathbb{Z}\times\left\{ R(i)_{2}+k\left\lfloor \sqrt{\lambda m}\right\rfloor \right\} \\
S_{k} & :=\mathbb{Z}\times\left(R(i)_{2}+\left](k-1)\left\lfloor \sqrt{\lambda m}\right\rfloor ,(k+1)\left\lfloor \sqrt{\lambda m}\right\rfloor \right[\:\right)\\
H_{k} & :=\{ v\in S_{k}:|v_{2}|\geq|R(i)_{2}|\}.\end{align*}
Denote by $T_{k}$ the first time (after $i$) when $R$ hits $L_{k}$
(for the purpose of the definition of $T_{k}$ we extend the walk
to infinity). Next denote by $T_{k}^{*}>T_{k}$ the first time after
$T_{k}$ when $R$ exits the strip $S_{k}$. Examine now the event
$G_{k}$ that $R$ has $\geq c(\lambda m)^{3/4}$ relative tan points
in the half strip $H_{k}$ in the time interval $[T_{k},T_{k}^{*}]$.
Translation and reflection symmetry shows that the $G_{k}$ are i.i.d.
It is known that $\mathbb{P}(G_{k})\geq c_{2}$, see \cite{BW03}.
Further, if $\lambda$ is sufficiently small then $\mathbb{P}(T_{k}^{*}-T_{k}>m)\leq\frac{1}{2}c_{2}$
uniformly in $m>1/\lambda$. Fix $\lambda$ to satisfy this requirement
and denote $G_{k}^{*}=G_{k}\cap\{ T_{k}^{*}-T_{k}\leq m\}$ so that
$\mathbb{P}(G_{k}^{*})\geq c$.

The events $G_{k}^{*}$ are also independent so if we denote by $B_{3}^{\pm}$
the event \[
\neg G_{k}^{*}\textrm{ for all }k\in\pm\left]\frac{\log^{2}n}{\sqrt{\lambda}},2\frac{\log^{2}n}{\sqrt{\lambda}}-1\right[\]
then $\mathbb{P}(B_{3}^{\pm})\leq Ce^{-c\log^{2}n}$ (here we consider
$\lambda$ as a constant and allow $C$ and $c$ to depend on it).
The lemma is now finished since if none of the four bad events $B_{1}$,
$B_{2}$, $B_{3}^{+}$, $B_{3}^{-}$ happened the claim holds. Indeed,
$\neg B_{2}$ implies that either \[
T_{\left\lfloor 2\lambda^{-1/2}\log^{2}n\right\rfloor }<i+m\log^{6}n-m\quad\textrm{or}\quad T_{-\left\lfloor 2\lambda^{-1/2}\log^{2}n\right\rfloor }<i+m\log^{6}n-m.\]
 Assume that the first happened. Then if $B_{3}^{+}$ did not happen
then some $G_{k}^{*}$ happened and by the definition of $G_{k}^{*}$
we can denote $j:=T_{k}^{*}$ and get what we want ($\neg B_{1}$
is used to show $j\geq i+m$). Hence with probability $\geq1-Ce^{-c\log^{2}n}$
we found a $j$ for our $i$. Summing over $i$ we are done.
\end{proof}
\begin{lem}
\label{lem:SRW2}With probability $\geq1-Ce^{-c\log^{2}n}$ one has
that for any $i\neq j$, $|R(j)_{1}-R(i)_{1}|\leq\log n\sqrt{i-j}$.
\end{lem}
This follows immediately from the Chernoff bound.

\begin{lem}
\label{lem:Em3/4}Let $E$ be an excited random walk of length $2n$
starting from some $\mathcal{V}\subset\left]-\infty,-n^{5/8}\right]\times\mathbb{Z}$
and some vertex $E(0)\in\mathbb{Z}^{2}$. Let $m\geq n^{15/16}$.
Then with probability $\geq1-C_{1}e^{-c\log^{2}n}$ one has that either
\begin{enumerate}
\item \label{enu:Rn<0}$E(n)_{1}<0$; or
\item \label{enu:Em3/4}For every $n\leq i\leq2n-m\log^{6}2n$, \[
E\big(i+\left\lfloor m\log^{6}2n\right\rfloor \big)_{1}-\max_{j\leq i}E(j)_{1}\geq cm^{3/4}.\]

\end{enumerate}
\end{lem}
\begin{proof}
We may assume $n$ is sufficiently large (for $n$ small choosing
$C_{1}$ sufficiently large will render the lemma true trivially).
If $m>n\log^{-6}2n$ then \ref{enu:Em3/4} holds vacuously so assume
the opposite. Couple $E$ to a simple random walk $R$ as in the introduction.
Let $B_{1}$ and $B_{2}$ be the bad events of lemmas \ref{lem:imj}
and \ref{lem:SRW2} for the walk $R$, i.e.~$B_{1}$ is the event
that for some $i\in\{0,\dotsc,\left\lfloor 2n-m\log^{6}2n\right\rfloor \}$
there aren't enough relative tan points and $B_{2}$ is the event
that for some $i\neq j$ $|R(i)_{1}-R(j)_{1}|$ is very large. Let
$T_{k}$ be the $k$-th time when $E$ reached a new vertex and let
$\xi_{k}$ be $E(T_{k+1})-E(T_{k})-(R(T_{k+1})-R(T_{k}))$. $\xi_{k}$
is a vector but since it can take only the values $(0,0)$ and $(2,0)$
we will consider it as a scalar. Let $B_{3}$ be the event that for
some $k\leq2n$, \[
\sum_{l=k}^{k+\left\lfloor c_{1}m^{3/4}\right\rfloor }\xi_{l}\leq c_{3}m^{3/4}.\]
Now, the $\xi_{k}$-s are independent hence it is easy to see that
for $c_{3}$ sufficiently small for any $k$ the probability for this
is $<Ce^{-cm^{3/4}}$ and then so will be their sum. In other words,
$\mathbb{P}(B_{3})<Ce^{-cm^{3/4}}$.

The lemma would be proved if we show that $\neg(B_{1}\cup B_{2}\cup B_{3}\cup\textrm{\ref{enu:Rn<0}})\Rightarrow\textrm{\ref{enu:Em3/4}}$.
This will be done inductively, and the first step is to choose a number
$1\leq i_{1}\leq n$ as follows. We divide into two cases according
to whether $E(0)_{1}>-\frac{1}{2}n^{5/8}$ or not. If $E(0)_{1}>-\frac{1}{2}n^{5/8}$
we choose $i_{1}=0$. Otherwise we note that if $E(0)_{1}\leq-\frac{1}{2}n^{5/8}$
and $E(n)_{1}>0$ then we can choose some $i_{1}\in[0,\dotsc,n]$
with the properties that $E\left(i_{1}\right)_{1}>-\frac{1}{2}n^{5/8}$
and $E\left(i_{1}+m\right)_{1}-E\left(i_{1}\right)_{1}>\frac{1}{4}mn^{-3/8}\geq\frac{1}{4}n^{9/16}$
for $n$ sufficiently large.

Now, by $\neg B_{2}$ we see that for all $j\leq i_{1}$ one has $R(j)_{1}\leq R\left(i_{1}\right)_{1}+\log2n\sqrt{2n}$
and the coupling implies the same for $E$. Similarly we get that
for all $j\geq i_{1}+m$, $E(j)_{1}\geq E\left(i_{1}+m\right)_{1}-\log2n\sqrt{2n}$.
Hence we see that (for $n$ sufficiently large) \[
(\mathcal{V}\cup E[1,i_{1}])\cap E[i_{1}+m,2n]=\emptyset.\]
 This means that any tan point of $R$ relative to $i_{1}$ after
$i_{1}+m$ is a point where $E$ encountered a new vertex. By $\neg B_{1}$
we know that for some $i_{2}\in[i_{1}+m,i_{1}+m\log^{6}n-m]$ we would
have $\geq c_{1}m^{3/4}$ relative tan points in $[i_{2},i_{2}+m]$
and hence $\geq c_{1}m^{3/4}$ new points for $E$. By $\neg B_{3}$
we see that \[
E\left(i_{2}+m\right)_{1}-E\left(i_{2}\right)_{1}>c_{3}m^{3/4}+R\left(i_{2}+m\right)_{1}-R\left(i_{2}\right)_{1}\]
and by $\neg B_{2}$ \[
\geq c_{3}m^{3/4}-\log2n\sqrt{2n}>cm^{3/4}>{\textstyle \frac{1}{4}}n^{9/16}\]
 for $n$ sufficiently large.

Now we can repeat the argument of the last paragraph with $i_{1}$
replaced by $i_{2}$. We get a sequence of $i$-s satisfying (except
possibly $i_{1}$)\[
E(i_{j}+m)-E(i_{j})>cm^{3/4}\]
and $i_{j}\leq i_{j-1}+m\log^{6}n-m$. This implies (again with $\neg B_{2}$)
\ref{enu:Em3/4} and the lemma.
\end{proof}
\begin{lem}
\label{lem:B}Let $B_{1},\dotsc,B_{n}$ be i.i.d.~Bernoulli (i.e.~$0$-$1$)
variables with $\mathbb{P}(B_{i}=1)=\epsilon$. Then $\mathbb{P}(\sum B_{i}>k)\leq2(n\epsilon)^{k}$.
\end{lem}
This is a straightforward calculation (and a very rough estimate to
boot --- we will only use it for $\epsilon\ll1/n$ where it is rather
close to the truth).

\begin{lem}
\label{lem:induct}Let $E,n,\mathcal{V}$ and $m$ be as in lemma
\ref{lem:Em3/4} and assume in addition that $m\geq2n^{15/16}\log^{6}2n$.
Assume also that one knows that an excited random walk of length $2m$
starting from any $\mathcal{W}\subset\left]-\infty,-m^{5/8}\right]\times\mathbb{Z}$
and any vertex in $\mathbb{Z}^{2}$ satisfies, with probability $\geq1-\epsilon$,
that either \begin{align*}
E(m)_{1} & <0\textrm{; or}\\
E(2m)_{1}-E(m)_{1} & >\mu.\end{align*}
(here end the assumptions of the lemma). Then with probability $\geq1-Ce^{-c\log^{2}n}-4(n\epsilon)^{\left\lfloor \log n\right\rfloor }$,
either \begin{align}
E(n)_{1} & <0\textrm{; or}\label{eq:En1<0}\\
E(2n)_{1}-E(n)_{1} & >\mu(\left\lfloor n/m\right\rfloor -2\log n-1).\label{eq:E2nEn>}\end{align}

\end{lem}
\begin{proof}
We may assume $n$ is sufficiently large. Let $k\in\{0,\dotsc,\left\lfloor n/m\right\rfloor -2\}$
and let $B^{k}$ be the event that \begin{gather*}
E(n+(k+1)m)_{1}>\max\left\{ \max_{j\leq n+km}E(j)_{1},\max_{v\in\mathcal{V}}v_{1}\right\} +m^{5/8}\\
E(n+(k+2)m)_{1}-E(n+(k+1)m)\leq\mu.\end{gather*}
 We translate by $-E(n+(k+1)m)$ and get from the assumption that
\[
\mathbb{P}(B^{k}\,|\, E[0,\dotsc,n+km])\leq\epsilon.\]
Since $B^{0},\dotsc,B^{k-2}$ depend only on $E[0,\dotsc,km]$ then
we get that the even $B^{k}$-s are stochastically dominated by a
sequence of i.i.d.~$\epsilon$-Bernoulli variables. The odd $B^{k}$
satisfy the same. Hence by lemma \ref{lem:B}, \[
\mathbb{P}(\# B^{k}>2\log n)\leq4(\left\lceil \left\lfloor n/m\right\rfloor /2\right\rceil \epsilon)^{\left\lfloor \log n\right\rfloor }<4(n\epsilon)^{\left\lfloor \log n\right\rfloor }.\]
Denote this event by $B_{1}$.

We now apply lemma \ref{lem:Em3/4} with $m_{\textrm{lemma \ref{lem:Em3/4}}}=l=\left\lfloor m\log^{-6}n\right\rfloor $
(and the same $n$). If $n$ is sufficiently large then $l\geq n^{15/16}$
and lemma \ref{lem:Em3/4} may indeed be applied. We get that with
probability $>1-Ce^{-c\log^{2}n}$ we have either (\ref{eq:En1<0})
or for every $n\leq i\leq2n-l\log^{6}2n$, \begin{equation}
E\left(i+\left\lfloor l\log^{6}2n\right\rfloor \right)_{1}-\max_{j\leq i}E(j)_{1}\geq cl^{3/4}.\label{eq:llog62n}\end{equation}
Since $\left\lfloor l\log^{6}2n\right\rfloor \simeq m$ (the difference
is $\leq\log^{6}2n+1$) we can replace the first by the second in
(\ref{eq:llog62n}) and pay only in the constant. Denote therefore
by $B_{2}$ the event that $\neg(\ref{eq:En1<0})$ and also \[
\exists i\in\{ n,\dotsc,2n-m\}:E(i+m)_{1}-\max_{j\leq i}E(j)_{1}<c_{4}l^{3/4}.\]
Then for $c_{4}$ sufficiently small we have $\mathbb{P}(B_{2})\leq Ce^{-c\log^{2}n}$.

We claim that $\neg(B_{1}\cup B_{2}\cup(\textrm{\ref{eq:En1<0}}))\Rightarrow(\textrm{\ref{eq:E2nEn>}})$,
which will finish the lemma. This however, is clear: $E(n)_{1}\geq0$
and $\neg B_{2}$ give that for all $k\in\{0,\dotsc,\left\lfloor n/m\right\rfloor -2\}$,
\[
E(n+(k+1)m)_{1}\geq\max\left\{ \max_{j\leq n+km}E(j)_{1},\max_{v\in\mathcal{V}}v_{1}\right\} +cl^{3/4}\]
 so for $n$ sufficiently large (so that $cl^{3/4}>m^{5/8}$) $\neg B^{k}$
will give \[
E(n+(k+2)m)_{1}>E(n+(k+1)m)+\mu.\]
 Only $\left\lfloor 2\log n\right\rfloor $ $B^{k}$-s are allowed
to fail, so \[
\sum_{k:\neg B^{k}}E(n+(k+2)m)_{1}-E(n+(k+1)m)_{1}>\mu(\left\lfloor n/m\right\rfloor -2\log n-1).\]
While if $B^{k}$ does not hold we can still use $\neg B_{2}$ to
get \[
E(n+(k+2)m)_{1}-E(n+(k+1)m)_{1}>0\]
and we are done.
\end{proof}
\begin{lem}
\label{lem:main}With the notations of lemma \ref{lem:Em3/4}, with
probability $\geq1-Ce^{-c\log^{2}n}$, either 
\begin{enumerate}
\item $E(n)_{1}<0$; or
\item $E(2n)_{1}-E(n)_{1}>cn$.
\end{enumerate}
\end{lem}
\begin{proof}
This follows by an inductive application of lemma \ref{lem:induct},
but one has to be careful with the parameters. Let therefore $K$
and $k$ be two parameters which will be fixed later. Define $\alpha_{n}=\alpha_{n}(K,k)$
as the maximal number such that\[
\mathbb{P}(\{ E(n)_{1}>0\}\cap\{ E(2n)_{1}-E(n)_{1}\leq\alpha_{n}n\})\leq Ke^{-k\log^{2}n}\quad\forall\mathcal{V},E(0)\]
(as usual we assume $\mathcal{V}\subset\left]-\infty,-n^{5/8}\right]\times\mathbb{Z}$).
We wish to estimate $\alpha_{n}$. First we check what lemma \ref{lem:Em3/4}
has to say about $\alpha_{n}$. Choosing $m=\left\lfloor n\log^{-6}2n\right\rfloor $
(which can be done if $n$ is sufficiently large) and $i=n$ we get
that with probability $>1-Ce^{-c\log^{2}n}$, either $E(n)_{1}<0$
or\begin{align*}
E(2n)_{1}-E(n)_{1} & \geq E\left(n+\left\lfloor m\log^{6}2n\right\rfloor \right)_{1}-E(n)_{1}-\log^{6}2n-1\geq\\
 & \geq cm^{3/4}-\log^{6}2n-1>cn^{5/8}.\end{align*}
 As usual we can remove the assumption that $n$ is sufficiently large
and pay only in the constants. In other words, if $K$ is sufficiently
large and $k$ is sufficiently small then $\alpha_{n}\geq c(K,k)n^{-3/8}$
for all $n$. 

Next we translate lemma \ref{lem:induct} to $\alpha_{n}$ notations
and it now goes: if we have \begin{equation}
Ce^{-c\log^{2}n}+4\left(nKe^{-k\log^{2}m}\right)^{\left\lfloor \log n\right\rfloor }\leq Ke^{-k\log^{2}n}\label{eq:KkCc}\end{equation}
then \[
\alpha_{n}\geq\alpha_{m}\left(1-\frac{(2\log n-2)m}{n}\right).\]
 It is easy to see that for $K$ sufficiently large and $k$ sufficiently
small (\ref{eq:KkCc}) will hold for $n$ sufficiently large (all
bounds depend only on the $C$ and $c$ in (\ref{eq:KkCc})). Fix
$K$ and $k$ to satisfy both requirements. We get that for $n$ sufficiently
large, \begin{equation}
\alpha_{n}\geq\alpha_{m}\left(1-n^{-1/32}\right)\quad m=\left\lceil 2n^{15/16}\log^{6}2n\right\rceil .\label{eq:alpha n m}\end{equation}
Let $N$ satisfy that for all $n>N$ (\ref{eq:alpha n m}) holds and
in addition $m<\frac{1}{2}n$. We easily get\[
\alpha_{n}>c\min_{m\leq N}\alpha_{m}>cN^{-3/8}\quad\forall n.\]
and we are done.
\end{proof}

\begin{proof}
[Proof of the theorem]This is an immediate corollary from lemma \ref{lem:main}.
Take $\mathcal{V}=\emptyset$ and $E(0)=(n+1,0)$ (so that $E(n)_{1}>0$
always) and then translate by $-n-1$. We get with probability $>1-Ce^{-c\log^{2}n}$,
\[
E(2n)>E(n)+cn\geq R(n)+cn\]
 and since with probability $>1-Ce^{-c\log^{2}n}$ we have $R(n)>-\log n\sqrt{n}$
the theorem is proved.
\end{proof}

\end{document}